\documentclass{amsart}

\usepackage{amsthm,amsmath,amsfonts}
\usepackage[all]{xy}
\usepackage{psfrag}
\usepackage[dvips]{graphicx}
\usepackage{epsf, subfigure, verbatim}
\usepackage[latin1]{inputenc}
\usepackage[T1]{fontenc}
\usepackage{pslatex}

\makeatletter
\def\th@alexnormal{%
\let\thm@indent\noindent 
\thm@headfont{\bfseries}
\normalfont
}
\makeatother

\makeatletter
\def\th@alexit{%
\let\thm@indent\noindent 
\thm@headfont{\bfseries}
\normalfont
\fontshape{it}
\selectfont
}
\makeatother

\theoremstyle{alexit}
\newtheorem{theorem}{Theorem}[section]
\newtheorem{proposition}[theorem]{Proposition}
\newtheorem{corollary}[theorem]{Corollary}
\newtheorem{lemma}[theorem]{Lemma}

\theoremstyle{alexnormal}
\newtheorem{definition}[theorem]{Definition}
\newtheorem{remark}[theorem]{Remark}

\begin{document}
\title{Relative homological linking in critical point theory}
\author[A. Girouard]{Alexandre Girouard}

\address{D\'epartement de Math\'ematiques et
Statistique, Universit\'e de Montr\'eal, C. P. 6128,
Succ. Centre-ville, Montr\'eal, Canada H3C 3J7}
\email{girouard@dms.umontreal.ca}

\begin{abstract}
  A relative homological linking of pairs is proposed. It is shown to
  imply homotopical linking, as well as earlier non-relative notion of
  homological linkings. Using Morse theory we prove a simple
  ``homological linking principle'', thereby generalizing and
  simplifying many well known results in critical point theory.
\end{abstract}

\maketitle

\section*{Introduction}
The use of linking methods in critical point theory is rather new. It
was implicitely present in the work of Ambrosetti and Rabinowitz
\cite{ar:1} in the early 70's as well as in the work of Benci and Rabinowitz \cite{benci:1}.
The first explicit definition was given by Ni in 1980 \cite{ni:1}.

\begin{definition}[Classical Homotopical Linking]
  Let $A\subset B$ and $Q$ be subspaces of a topological space $X$ 
  such that the pair $(B,A)$ is homeomorphic to $(D^n,S^{n-1})$. Then 
  \emph{$A$ homotopically links $Q$} if for each deformation
  $\eta:[0,1]\times B\rightarrow X$ fixing $A$, $\eta(1,B)\cap Q \neq
  \emptyset$.
\end{definition}

In the early 80's, homological linking was introduced in critical point theory
(see Fadell \cite{fadell:1}, Benci \cite{benci:2} and Chang \cite{chang:2} for instance).

\begin{definition}[Classical Homological Linking]
  Let $A$ and $S$ be non-empty disjoint subspaces in a topological space $X$.
  Then \emph{$A$ homologically links $S$} if the inclusion of $A$ in $X\setminus S$ induces a
  non-trivial homomorphism in reduced homology.
\end{definition}

In her 1999's article \cite{frigon:1}, Frigon generalized homotopical linking to pairs of subspaces.
\begin{definition}[Relative Homotopical Linking]
  Let $(B,A)$ and $(Q,P)$ be two pairs of subspaces in a topological space $X$
  such that $B\cap P=\emptyset$ and  $A\cap Q=\emptyset$. Then 
  \emph{$(B,A)$ homotopically links $(Q,P)$} 
  if for each deformation $\eta:[0,1]\times B\rightarrow
  X$ fixing $A$ pointwise, $\eta(1,B)\cap Q=\emptyset\Rightarrow \exists t\in
  ]0,1], \eta(t,B)\cap P\neq \emptyset$.
\end{definition}
The classical definition corresponds to the case where
$(B,A)\cong (D^n,S^{n-1})$ and $P=\emptyset$.

The goal of this article is to propose a similar generalization for
homological linking.
In section 1.1 we explore the properties of this new homological
linking and in 1.2 we give some detailed examples.
In section 2 we interpret homotopical linking as an obstruction to
factoring certain homotopy through homotopically trivial pairs. It
becomes clear from this point of view that homological linking is
stronger than homotopical linking. 
Our definition of homological linking fits very nicely with Morse theory. We exploit this
in section 3 to derive a new linking principle (see
\ref{homo_link_prin}) for detecting and locating critical
points. Despite its simplicity, the idea is quite fruitful. Close
analog to the Mountain Pass Theorem of Ambrosetti and Rabinowitz
\cite{ar:1} as well as to the Saddle Point Theorem of Rabinowitz
\cite{rabinowitz:1} are easy corollaries. In Proposition
\ref{FrigHomotop}, we also obtain a homological version  of the
generalized saddle point theorem of Frigon \cite{frigon:1}. In section
4, some multiplicity results are studied.

Our approach has many advantages: each critical point is detected by a different linking,
stability type is directly available (i.e. critical groups are known) and last but not least, 
the proofs are easy. However, it also has a disadvantage: working with Morse theory requires
more regularity  than using  a ``min-max'' method for example.
It might appear as if the content of this paper is extremely easy. We agree
with this point of view. In fact, it is rather surprising to see that so many of the
classical results of critical point theory are straightforward consequences of this new
definition of homological linking. 

This paper is an extension of the author's master's thesis \cite{g:1}. He
would like to express his most sincere thanks to his advisor, Marlène
Frigon.

\tableofcontents

\section{Homological linking}

\subsection{Definition and properties}
The principal contribution of this article is the following definition.
\begin{definition}[Relative Homological Linking]
  Let $(B,A)$ and $(Q,P)$ be pairs of subspaces in a topological
  space $X$.
  Then {\em $(B,A)$ homologically links $(Q,P)$  in $X$}
  if $(B,A) \subset (X\setminus P, X\setminus Q)$
  and if this inclusion induces
  a non-trivial homomorphism in reduced homology.
  Given integers $q, \beta \geq 0$, we say that
  {\em $$(B,A)\ (q,\beta)\mbox{-links } (Q,P)\mbox{ in } X$$} if 
  the above inclusion induces a homomorphism of rank~$\beta$ on the
  $q$-th reduced homology groups.
\end{definition}


\begin{remark}
  For notational convenience, a topological pair $(B,\emptyset)$ will be
  identified with the space $B$.
\end{remark}

\begin{remark}
  The classical definition corresponds to the case
  $A\ (q,\beta)\mbox{-links } (X,Q)$ and $\beta>0$.
\end{remark}

\begin{remark}\label{betti}
  For any space $X$, $X\ (q,b_q(X))\mbox{-links } X$ in $X$, where $b_q(X)$ is the $q$-th reduced Betti number of $X$.
  Thus our linking contains as much information as Betti numbers.
\end{remark}

The next proposition and it's corollary shows that in many situations,
it suffices to consider linking locally to deduce a global linking situation.

\begin{proposition}
  Let $\mathcal{O}$ be an open subset of $X$. If\\ $A,B,P,Q\subset
  \mathcal{O}$ with $Q$ closed, then
  \begin{gather*}
    (B,A)\ (q,\beta)\mbox{-links } (Q,P)\mbox{ in } X\\
    \Leftrightarrow\\
    (B,A)\ (q,\beta)\mbox{-links } (Q,P)\mbox{ in } \mathcal{O}.
  \end{gather*}
\end{proposition}
\begin{proof}
  Since $\mathcal{O}^c$ is closed and $X\setminus Q$ is open
  in $X\setminus P$, the excision axiom applies to 
  $$\mathcal{O}^c \subset X\setminus Q \subset X\setminus P.$$
  It follows that the the bottom line of the following commutative diagram is an isomorphism.
  $$\xymatrix{
    \tilde{H}_q(B,A)\ar[d]^i\ar[rd]^j & \\
    \tilde{H}_q(\mathcal{O}\setminus P, \mathcal{O}\setminus Q)\ar[r]^{\cong} &
    \tilde{H}_q(X\setminus P, X\setminus Q)
    }$$
  Hence, $\mbox{rank } j=\mbox{rank } i$.
\end{proof}

\begin{corollary}
  Let $\mathcal{O}$ be the domain of a chart on a manifold $M$. If the
  pair $(B,A)$ links the pair $(Q,P)$ in $\mathcal{O}$, with $Q$
  closed, then $(B,A)$ also links the pair $(Q,P)$ in $M$.
\end{corollary}

The two following theorems show how some simple linking situations
lead to new linkings.

\begin{theorem}\label{thm_linking_1}
  If $A\ (q,\beta)\mbox{-links } (X,Q)$ and $A\ (q, \delta)\mbox{-links } (X,X\setminus B)$ in $X$
  for some $\delta<\beta$ then $(B,A)\ (q+1,\mu)\mbox{-links } Q$ in $X$ for some $\mu\geq\beta-\delta$.
\end{theorem}
\begin{proof}
  It follows from the commutativity of
  $$\xymatrix{
    \tilde{H}_{q+1}(B,A) \ar[r]^{\Delta_1}\ar[d]^\alpha&\tilde{H}_q(A)\ar[r]^k\ar[d]^i&\tilde{H}_q(B)\ar[d]\\
    \tilde{H}_{q+1}(X,X\setminus Q) \ar[r]^{\Delta_2}&\tilde{H}_q(X\setminus Q)\ar[r]&\tilde{H}_q(X)
    }$$
  that
  \begin{align*}
    \mu:=\mbox{rank }\alpha &\geq\mbox{rank } \Delta_2\circ  \alpha =\mbox{rank }i\circ \Delta_1\\
    &\geq \mbox{rank }\Delta_1 - \dim(\ker i)\\
    &= \mbox{rang }\Delta_1 - (\dim \tilde{H}_q(A) - \mbox{rank }i)\\
    &= \mbox{rank }i+\mbox{rank }\Delta_1 -\dim \tilde{H}_q(A)\\
    &= \mbox{rank }i+\mbox{rank }\Delta_1 -(\mbox{rank }k + \dim(\ker k)).
  \end{align*}
  By exactness, $\mbox{rank }\Delta_1 = \dim(\ker k)$, thus
  $$\mu \geq \mbox{rank }i - \mbox{rank }k = \beta - \delta.$$
\end{proof}

 \begin{theorem}\label{thm_linking_2}
   If $B\ (q,\beta)\mbox{-links } (X,P)$ and 
   $X\setminus Q\ (q,\delta)\mbox{-links } (X,P)$ for some
   $\delta<\beta$, then $B\ (q,\mu)\mbox{-links }(Q,P)$ in $X$
   for some $\mu\geq \beta-\delta$.
 \end{theorem}
\begin{proof}
From the commutativity of
$$\xymatrix{
  &\tilde{H}_q(B)\ar[r]^\cong\ar[d]^i&\tilde{H}_q(B,\emptyset)\ar[d]^\alpha\\
  \tilde{H}_q(X\setminus Q)\ar[r]^k&\tilde{H}_q(X\setminus P)\ar[r]^j&\tilde{H}_q(X\setminus P,X\setminus Q)
  }$$
it follows that
\begin{align*}
  \mu=\mbox{rank }\alpha &= \mbox{rank }j\circ i\\
  &\geq \mbox{rank }i -\dim(\ker j)\\
  &=\mbox{rank }i - \mbox{rank }k\\
  &= \beta - \delta.
\end{align*}
\end{proof}

\subsection{Examples of linking}
Our definition permits to obtain new situations of linking and to recover others already known.
In particular, in Propositions \ref{enlacement_1}, \ref{enlacement_2} and \ref{enlacement_3} 
we present linking situations equivalent to those already studied by
Perera in \cite{perera:1} using a non relative definition of homological
linking.

Let $E$ be a Banach space. Given a direct sum decomposition
$E = E_1 \oplus E_2$, $B_i$ denotes the closed ball in $E_i$ and $S_i$
its relative boundary ($i=1,2$).

\begin{proposition}\label{enlacement_1}
  Let $e\in E$, $\|e\|>1$.
  Then $\{0, e\}\ (0,1)\mbox{-links }(E,S)$ in $E$.
\end{proposition}
\begin{proof}
 The map $r:E\setminus S \rightarrow \{0,e\}$ defined
 by
 $$r(x)=\left\{
 \begin{array}{ll}
   0&\mbox{if } \|x\| < 1,\\
   e&\mbox{if } \|x\| > 1.
 \end{array}
 \right.
 $$
 is a retraction. That is, the following diagram commutes
 $$\xymatrix{
   E\setminus S \ar[r]^r&\{0,e\}\\
   \{0,e\}\ar[ru]_{\mbox{id}}\ar[u]&
   }$$
 It follows that the inclusion of $\{0,e\}$ in $E\setminus S$ is of
 rank 1 in reduced homology.
\end{proof}

\begin{proposition}\label{enlacement_2}
  Let $E=E_1\oplus E_2$ with $k=\dim E_1 \in\ ]0, \infty[$.
  Then $$S_1 \ (k-1,1)\mbox{-links } (E,E_2)$$ in $E$.
\end{proposition}

\begin{proof}
 The long exact sequence induced by $S_1~\subset~E\setminus~E_2$ is
 $$\cdots\rightarrow \tilde{H}_k(E\setminus E_2, S_1) \rightarrow
 \tilde{H}_{k-1}(S_1) \stackrel{i}{\rightarrow}
 \tilde{H}_{k-1}(E\setminus E_2)\rightarrow \cdots$$
 Because $E\setminus E_2$ strongly retract on $S_1$,
 $H_k(E\setminus E_2,S_1)=0$. It follows that
 $\mbox{rank }i=\dim\tilde{H}_{k-1}(S_1)=1$.
\end{proof}

\begin{proposition}\label{enlacement_3}
  Let $E=E_1\oplus E_2$ with $k=\dim E_1 \in\ ]0, \infty[$
  and let $e \in E_2$ be of unit length. 
  Let $A=\partial (B_1\oplus [0,2]e)$ in $E_1\oplus \mathbb{R}e$.
  Then $A \ (k,1)\mbox{-links } (E,S_2)$ in $E$.
\end{proposition}

\begin{proof}
 Let $P:E\rightarrow E_1$ be the projection on $E_1$ and
 $r:E\setminus S_2 \rightarrow (E_1\oplus \mathbb{R}e) \setminus \{e\}$
 be defined by $r(x)=P(x)+\|x-P(x)\|e$. 
 Let's make sure $\{e\}$ really is omitted by $r$.
 Suppose $x\in E$ is such that $P(x)+\|x-P(x)\|e=e$. Then $P(x)=0$ and
 $1=\|x-P(x)\|=\|x\|$. In other words,  $x\in E_2$ and $\|x\|=1$ wich
 is impossible for $x$ in the domain of $r$.
 Let $i$ be the inclusion of $A$ in $E\setminus S_2$.
 If $i_k:\tilde{H}_k(A)\rightarrow \tilde{H}_k(E\setminus S_2)$ is
 null, then so is 
 $$r_k\circ i_k:\tilde{H}_k(A)\rightarrow \tilde{H}_k((E_1\oplus
 \mathbb{R}e)\setminus \{e\}).$$
 However, $r\circ i$ is the inclusion of $A$ in $(E_1\oplus
 \mathbb{R}e)\setminus \{e\}$ and
 $(E_1\oplus \mathbb{R}e)\setminus \{e\}$ strongly retract on 
 $A$. Thus $\tilde{H}_*((E_1\oplus \mathbb{R}e)\setminus
 \{e\},A)\cong 0$. It then follows from the long exact sequence
 induced by the inclusion $r\circ i$ of $A$ in $(E_1\oplus \mathbb{R}e)\setminus \{e\}$
 $$0=\tilde{H}_{k+1}((E_1\oplus \mathbb{R}e)\setminus \{e\},A)\rightarrow 
 \tilde{H}_k(A)\overset{{r_k\circ i}_k}{\rightarrow}
 \tilde{H}_k(E_1\oplus \mathbb{R}e\setminus \{e\})$$
 that ${r_k\circ i}_k$ is not trivial because $\tilde{H}_k(A)\cong
 \mathbb{K}$. Consequently $A\ (k,1)\mbox{-links } (E,S_2)$ in $E$, as
 was to be proved.
\end{proof}

Theorem \ref{thm_linking_1} and the previous linking situations give
rise to other linkings which are in fact the classical situations
treated in the litterature. Observe that, in these classical
situations, the pair $(Q,P)$ is always of the form $(Q,\emptyset)$ and
the pair $(B,A)$ always has $A\neq\emptyset$.

\begin{corollary}\label{cor_enlacement_1}
  Let $e\in E$ with $\|e\|>1$.
  Then $([0,e],\{0, e\})\ (1,1)$-links $S$ in $E$.
\end{corollary}

\begin{corollary}\label{cor_enlacement_2}
  Let $E= E_1 \oplus E_2$ with $k=\dim E_1 \in\ ]0, \infty[$.
  Then $$(B_1,S_1) \ (k,1)\mbox{-links } E_2$$ in $E$.
\end{corollary}

\begin{corollary}\label{cor_enlacement_3}
  Let $E= E_1 \oplus E_2$ with $k=\dim E_1 \in\ ]0, \infty[$ 
  and let $e \in E_2$ be of unit length.
  Let $B=B_1\oplus [0,2]e$ and $A=\partial B$ in
  $E_1\oplus \mathbb{R}e$.
  Then $(B,A)\ (k+1,1)\mbox{-links }S_2$ in $E$.
\end{corollary}
                                
By combining the linking situations of proposition
\ref{enlacement_1}, \ref{enlacement_2} and \ref{enlacement_3} with theorem
\ref{thm_linking_2}, we get a new familly of linking situations. 
These linking situation will be particularyly useful in applications to
critical point theory since they will allow us to relax the a priori estimates on $f$. 
For these linking, the pair $(B,A)$ is always of the form $(B,\emptyset)$ and the pair $(Q,P)$
always has $P\neq\emptyset$.

\begin{corollary}\label{cor_enlacement_4}
  Let $e\in E$, $\|e\|>1$.
  Then $$\{0, e\}\ (0,1)\mbox{-links }(B,S)$$
  in $E$.
\end{corollary}

\begin{corollary}\label{cor_enlacement_5}
  Let $E=E_1\oplus E_2$ with $k=\dim E_1 \in\ ]0, \infty[$ and let $e\in E_1$ be of unit
  length.
  Let $B=S_1, Q=E_2+[0,\infty[e$ and $P=E_2$.
  Then $$B \ (k-1,1)\mbox{-links } (Q,P)$$ in $E$.
\end{corollary}

\begin{corollary}\label{cor_enlacement_6}
  Let $E=E_1\oplus E_2$ with $k=\dim E_1 \in\ ]0, \infty[$
  and let $e \in E_2$ be of unit length. 
  Let $A=\partial (B_1\oplus [0,2]e)$ in $E_1\oplus \mathbb{R}e$.
  Then $A \ (k,1)\mbox{-links } (B_2,S_2)$ in $E$.
\end{corollary}


The two following propositions exhibit new homological linking
situations. From a homotopical point of view, they where studied by 
Frigon \cite{frigon:1}. These linking fully deserve to be called ``linking of pairs'' since for
both of them we have $A\neq\emptyset$ and $P\neq\emptyset$. A more geometrical argument is also
possible, but it is longuer.

\begin{proposition}\label{enlacement_4}
  Let $E=E_1\oplus E_2\oplus \mathbb{R}e$ with
  $e\in E$ of unit length and $k= \dim E_1 \in\ ]0,\infty[$. Let
  $B=B_1+e$, $A = S_1+e$, $Q=E_2+[0,\infty[e$ et $P=E_2$
  Then $(B,A) \ (k,1)\mbox{-links }(Q,P)$ in $E$.
\end{proposition}

\begin{proof}
  Let $\epsilon \in ]0,1[$ and 
  \begin{gather*}
    \hat{B}=B\cup (\epsilon B_1+]0,\infty[e+E_2),\\
    \hat{A}=\hat{B}\setminus (]0,\infty[e+E_2).
  \end{gather*}
  Since $B$ (resp. $A$) is a strong deformation retract of $\hat{B}$
  (resp. $\hat{A}$), the inclusion $(B,A)\rightarrow (\hat{B},\hat{A})$
  induces an isomorphism $H_k(B,A)\cong H_k(\hat{B},\hat{A})$. Let
  $$U=(E\setminus P)\setminus \hat{B} \subset E\setminus Q\subset
  E\setminus P,$$
  and observe that $\overline{U}\subset \mbox{int }(E\setminus Q)$ in
  $E\setminus P$, $\hat{B}=(E\setminus P)\setminus U$ and
  $\hat{A}=(E\setminus Q)\setminus U$. Hence, by excision, the inclusion
  $(\hat{B},\hat{A})\rightarrow (E\setminus P,E\setminus Q)$ induces an
  isomorphism $H_k(\hat{B},\hat{A})\cong H_k(E\setminus P,E\setminus
  Q)$. The result follows from $H_k(B,A)\cong \mathbb{K}$.
\end{proof}

A similar argument leads to the following proposition.

\begin{proposition}\label{enlacement_5}
  Let $E= E_1 \oplus E_2$ with $k=\dim E_1 \in\ ]0, \infty[$.\\Then
  $(B_1, S_1) \ (k,1)\mbox{-links } (B_2,S_2)$ in $E$.
\end{proposition}

\section{Homotopical consequences of homological linking}
Let $(B,A)$ and $(Q,P)$ be pairs of subspaces in a topological
space $X$ such that $B\cap P=\emptyset$ and $A\cap Q=\emptyset$.
The following lemma shows that relative homotopical linking is an
obstruction to extension factoring through a homotopically trivial pair.
\begin{lemma}
  The following statements are equivalent.
  \begin{enumerate}
  \item The pair $(B,A)$ homotopicaly links $(Q,P)$,
  \item There exists no homotopy
    $\eta:[0,1]\times (B,A)\rightarrow (X\setminus P, X\setminus Q)$
    such that
    $\eta=id$ on $\{0\}\times B\cup [0,1]\times A$ making the
    following diagram commutative
    $$\xymatrix{
      (B,A)\ar[r]^{\eta_1}\ar[rd]_{\eta_1}& (X\setminus P, X\setminus Q)\\
      &(X\setminus Q,X\setminus Q)\ar[u]
    }$$
  \end{enumerate}
\end{lemma}


\begin{corollary}
  Homological linking implies homotopical linking.
\end{corollary}

\begin{remark}
  To see that homotopical linking doesn't imply homological linking,
  it is sufficient to consider $X=B=Q$ to be a singleton and
  $A=P=\emptyset$.
\end{remark}

\section{Homological linking principle}
Let $H$ be a Hilbert space and let $f\in C^2(H,\mathbb{R})$. 
The following notation is standard. Given $c\in \mathbb{R}$,
$f_c=\{p\in H \bigr| f(p)\leq c\}$ is a level set of $f$,
$K(f)=\{p\in H \bigr| f'(p)=0\}$ is the critical set of $f$,
$K_c(f)=K(f)\cap f^{-1}(c)$.

Throughout this section, the following hypothesis
are assumed,
\begin{itemize}
\item[(H1)] the Palais-Smale condition for $f$ holds. That is, each
  sequence $(x_n)_{n\in\mathbb{N}}$ such that $(f(x_n))$ is bounded
  and $f'(x_n)\rightarrow 0$ admits a convergent subsequence,
\item[(H2)] the set $K(f)$ of critical point of $f$ is discrete.
\end{itemize}
In particular, $f(K)$ is discrete and for each bounded interval $I$,
$K\cap I$ is compact.

Under these assumptions, there is a suitable Morse theory which is
well behaved (see \cite{CPA:1} for instance). We shall use the
following standard notation. Given $p\in K_c(f)$,
$$C_q(f,p):=H_q(f_c,f_c\setminus\{p\})$$ 
is the $q$-th critical group of $f$ at $p$. Let $a<b$ be two regular values
of $f$, $$\mu_q(f_b,f_a):=\underset{p\in K(f)\cap f^{-1}[a,b]}{\sum} \dim C_q(f,p)$$
is the Morse number of the pair $(f_b,f_a)$.
The function $f$ is said to be a Morse function if its critical
points are all non-degenerate.

\begin{remark}
  Most of our results depend only on the Morse inequalities. 
  It is thus possible to use any other setting where they hold.
  For example, in \cite{corvellec:1} a Morse theory for continuous functions
  on metric spaces is presented. In applications to PDE, it may be necessary to use the
  Finsler structure approach of Chang \cite{chang:2} to apply the
  results in suitable Sobolev spaces. 
\end{remark}

The following theorem is an easy exercise and was probably first
observed by Marston Morse himself.
\begin{theorem}[homological linking principle]\label{homo_link_prin}
  Let $(B,A)$ and $(Q,P)$ be pairs of subspaces in $H$ and
  let $a<b$ be regular values of $f$ such that
  $(B,A) \subset (f_b,f_a) \subset (H\setminus P,H\setminus Q)$.
  If $(B,A)\ (q,\beta)$-links $(Q,P)$ in $H$ for some $\beta\geq 1$ then
  $f$ admits a critical point $p$ such that $a<f(p)<b$ and
  $C_q(f,p)\neq 0$. Moreover, if $f$ is a Morse function then it
  admits at least $\beta$ such points.
\end{theorem}
\begin{proof}
  It follows from commutativity of
  $$\xymatrix{
    \tilde{H}_q(B,A) \ar[r]\ar[d]&\tilde{H}_q(H\setminus P,H\setminus Q)\\
    \tilde{H}_q(f_b,f_a)\ar[ru]&
    }$$
  that $\dim \tilde{H}_q(f_b,f_a)\geq \beta$. Application of
  the weak Morse inequalities leads to \break$\mu_q(f_b,f_a)~\geq~\beta$ and to
  the first conclusion. The non-degeneracy condition leads to the
  second one.  
\end{proof}

\begin{remark}
  From Remark \ref{betti} and our linking principle we recover the weak Morse inequalities. This shows that our homological linking
  contains nearly as much information as classical Morse theory.
\end{remark}

\begin{lemma}
    Let $(B,A)$ and $(Q,P)$ be pairs of subspaces in $H$ such that
  \begin{gather*}
    \sup f(B) < \inf f(P),\\
    \sup f(A) \leq \inf f(Q).
  \end{gather*}
  If $(B,A)\ (q,\beta)$-links $(Q,P)$ in $H$ for some $\beta\geq 1$ then
  $\inf f(Q) \leq \sup f(B)$.
\end{lemma}
\begin{proof}
  Let the opposite be supposed: $\sup f(B)<\inf f(Q)$.  For each
  $n\in \mathbb{N}$, there exist regular values
  $a_n < b_n$ in $]\sup f(B), \sup f(B)+1/n[$.
  If $n$ is big enough, $\sup f(B)+1/n < \inf f(Q)\leq\inf f(P)$ so
  that $$(B,A)\subset (f_{b_n},f_{a_n}) \subset (X\setminus P, X\setminus Q).$$ 
  It follows from the homological linking principle that $f$ admits a
  critical value  $c_n\in\ ]a_n,b_n[$. The infinite sequence $(c_n)$
  converges to $c=\sup f(B)$ which must therefore be critical because
  the set of all critical values of $f$ is closed. This contradicts
  the fact that critical values must be isolated.
\end{proof}

The next theorem will be usefull for applications. In the next section, it will be used to
prove some multiplicity results.
\begin{theorem}\label{app_linking_principle}
  Let $(B,A)$ and $(Q,P)$ be pairs of subspaces in $H$ such that
  \begin{gather*}
    \sup f(B) < \inf f(P),\\
    \sup f(A) < \inf f(Q).
  \end{gather*}
  If $(B,A)\ (q,\beta)$-links $(Q,P)$ in $H$ for some $\beta\geq 1$ then
  $f$ admits a critical point $p$ such that
  $$\inf f(Q)\leq f(p)\leq \sup f(B)$$ and $C_q(f,p)\neq 0$.
  Moreover if $f$ is a Morse function then it admits at least $\beta$
  such points.
\end{theorem}

\begin{proof}
  By the preceding lemma,
  $$\sup f(A) < \inf f(Q) \leq \sup f(B) < \inf f(P).$$
  There exist regular values $a_n<b_n$ ($n\in \mathbb{N}$) 
  such that
  $$\sup f(A) < a_n < \inf f(Q) \leq \sup f(B) < b_n < \inf f(P)$$ and
  $a_n\rightarrow \inf f(Q)$, $b_n\rightarrow \sup f(B)$. By the
  linking principle, there must exist a sequence $(p_n)$ of critical
  points such that $C_q(f,p_n)\neq 0$ and such that the sequence
  $(c_n)=(f(p_n))$ satisfies $a_n<c_n<b_n$. Because critical values
  are isolated, $c_n\in [\inf f(Q), \sup f(B)]$ for $n$ big enough.
\end{proof}





The following result follows directly from Propositions~\ref{enlacement_5} and
Theorem~\ref{app_linking_principle}. As far as we know, this result is new.
\begin{theorem}\label{FrigHomotop}
  Let $H=H_1\oplus H_2$ with $k=\dim H_1<\infty$.
  If
  \begin{gather*}
    \sup f(S_1)<\inf f(B_2),\\
    \sup f(B_1)< \inf f(S_2)
  \end{gather*}
  then $f$ admits a critical point $p$ such that
  $$\inf f(S_2)\leq f(p) \leq \sup f(S_1)$$ and $C_k(f,p)\neq 0$.
\end{theorem}

\subsection{Multiplicity results}

By combining Corollaries \ref{cor_enlacement_3} and \ref{cor_enlacement_6} with Thorem
\ref{app_linking_principle}, we get a version of a well known multiplicity result
(see \cite{schechter:1} for instance). As before, we get extra information about the critical
groups.
\begin{proposition}
  Let $H=H_1\oplus H_2$ with $k=\dim H_1 \in\ ]0, \infty[$
  and $e \in H_2$ be of unit length.
  Let $B=B_1\oplus [0,2]e$ and $A=\partial B$ in in $H_1\oplus \mathbb{R}e$.
  If $f$ is bounded below on $B_2$ and if
  $$\sup f(A)<\inf f(S_2)$$ then $f$ admits two critical points
  $p_0\neq p_1$ such that 
  $$\inf(f(B_2)\leq f(p_0)\leq \sup f(A),$$
  $$\inf f(S_2)\leq f(p_1)\leq \sup f(B)$$
  and $C_k(f,p_0)\neq 0, C_{k+1}(f,p_1\neq 0)$. 
\end{proposition}
\begin{proof}
  Because
  \begin{gather*}
    \sup f(A) < \inf f(S_2)\\
    \sup f(\emptyset)=-\infty < \inf f(B_2)
  \end{gather*}
  and $A\ (k,1)\mbox{-links } (B_2,S_2)$, it follows from
  Theorem~\ref{app_linking_principle} that $f$ admits a critical point $p_0$
  such that $\inf f(B_2)\leq f(p_0) \leq \sup f(A)$ and
  $C_k(f,p_0)\neq 0$.
  Also, Corrolary \ref{cor_enlacement_3} says
  that $(B,A)\ (k+1,1)\mbox{-links } S_2$.
  Since 
  \begin{gather*}
    \sup f(B) < \infty=\inf f(\emptyset)\\
    \sup f(A) < \inf f(S_2)
  \end{gather*}
  it follows from Theorem~\ref{app_linking_principle} that $f$ admits 
  a critical point $p_1$ such that
  $\inf f(S_2)\leq f(p_1)\leq \sup f(B)$ and $C_{k+1}(f,p_1)\neq 0$.
  The inequality $$f(p_0)\leq \sup f(A) < \inf f(S_2)\leq f(p_1)$$
  insure that $p_0$ and $p_1$ are distinct.
\end{proof}

A similar argument using Corollaries \ref{cor_enlacement_2} and \ref{cor_enlacement_5}
leads to the next theorem. This result was already known to Perera \cite{perera:1}.
\begin{theorem}
  Let $H=H_1\oplus H_2$ with $k=\dim H_1 \in\ ]0, \infty[$ and let $e\in H_1$ be of unit
  length. If $f$ is bounded below on $H_1+[0,\infty[e$ and if 
  $$\sup f(S_1)<\inf f(H_2)$$ then $f$ admits two critical points $p_0\neq p_1$ such that
  $$\inf(f(H_1+[0,\infty[e))\leq f(p_0)\leq \max f(S_1),$$
  $$\inf f(H_2)\leq f(p_1)\leq \max f(B(0,1))$$
  and $C_{k-1}(f,p_0)\neq 0, C_k(f,p_1\neq 0)$.
\end{theorem}

\end{document}